\documentclass[12pt]{article}
\usepackage{amsmath}
\usepackage{amssymb}

\begin{document}
\title{Sharp Dimension Estimates of Holomorphic\\  Functions and Rigidity}
\author{Bing-Long Chen, Xiao-Yong Fu, Le Yin and Xi-Ping Zhu\\
{\small\it Department of Mathematics}\\
{\small\it Zhongshan University}\\
{\small\it Guangzhou, 510275, P.R.China}}
\date{}
\maketitle

\begin{abstract}\ \rm Let $M^n$ be a complete noncompact
K$\ddot{a}$hler manifold of complex dimension $n$ with nonnegative
holomorphic bisectional curvature. Denote by
$\mathcal{O}$$_d(M^n)$ the space of holomorphic functions of
polynomial growth of degree at most $d$ on $M^n$. In this paper we
prove that
$$dim_{\mathbb{C}}{\mathcal{O}}_d(M^n)\leq dim_{\mathbb{C}}{\mathcal{O}}_{[d]}(\mathbb{C}^n),$$
for all $d>0$, with equality for some positive integer $d$ if and
only if $M^n$ is holomorphically isometric to $\mathbb{C}^n$. We
also obtain sharp improved dimension estimates when its volume
growth is not maximal or its Ricci curvature is positive
somewhere.
\end{abstract}
\normalsize \vskip 1cm

\section{Introduction}\indent

In this paper we discuss the function theory of complete manifolds
with nonnegative curvature. Recall that the classical Liouville
theorem states that any bounded (or even just positive) harmonic
function is constant on Euclidean space. In [19], Yau extended the
classical Liouville theorem to complete noncompact Riemannian
manifolds with nonnegative Ricci curvature. It was further shown
by Cheng and Yau in [4] that any harmonic function with sublinear
growth on a complete noncompact Riemannian manifold with
nonnegative Ricci curvature must be constant. For the space
consisting of all harmonic functions of polynomial growth of
degree at most $d$ on a complete noncompact Riemanian manifold
$M^m$ of real dimension $m$ with nonnegative Ricci curvature
denoted by ${\mathcal{H}}_d(M^m)$, Yau asked in [20] if the
dimension of ${\mathcal{H}}_d(M^m)$ is finite for each positive
integer $d$, and if there holds
$$dim_{\mathbb{R}}{\mathcal{H}}_d(M^m)\leq
dim_{\mathbb{R}}{\mathcal{H}}_d({\mathbb{R}}^m)$$ for each
positive integer $d$, with equality for some positive integer $d$
if and only if $M^m$ is isometric to $\mathbb{R}^m$. The first
question was affirmativelly answered by Colding and Minicozzi [3],
and in [9] Li produced an elegant short proof. For the latter
question, the sharp upper bound estimate is still missing except
for the special cases $m=2$ or $d=1$ obtained by Li-Tam [11],[12],
and the rigidity part is only known for the special case $d=1$
obtained by Li [8] and Cheeger-Colding-Minicozzi [2].

In complex geometry, geometers studied (holomorphic) function
theory with the primary reason to prove a uniformization type
theorem in higher dimensions. In [16], Siu-Yau initiated a program
of using holomorphic functions of polynomial growth to
holomorphically embed complete K$\ddot{a}$hler manifold into
complex Euclidean spaces in an effort to generalize the classical
uniformization theorem. The well-known uniformization conjecture
for parabolic manifolds asks if every complete noncompact
K$\ddot{a}$hler manifold with nonnegative holomorphic bisectional
curvature is biholomorphic to a complex vector bundle over a
compact Hermitian symmetric space. For this sake, Yau [20]
proposed to study the spaces of holomorphic functions of
polynomial growth on complete noncompact K$\ddot{a}$hler manifolds
with nonnegative holomorphic bisectional curvature. Let us denote
by $\mathcal{O}_d(M^n)$ the space of holomorphic functions of
polynomial growth of degree at most $d$ on a complete noncompact
K$\ddot{a}$hler manifold $M^n$ of complex dimension $n$ with
nonnegative holomorphic bisectional curvature. The complex
counterpart of the above Yau's conjecture states if there holds
$$dim_{\mathbb{C}}{\mathcal{O}}_d(M^n)\leq
dim_{\mathbb{C}}{\mathcal{O}}_d({\mathbb{C}}^n)$$ for each
positive integer $d$, with equality for some positive integer $d$
if and only if $M^n$ is holomorphically isometric to
$\mathbb{C}^n$. In a recent beautiful work [14], Ni obtained the
following dimension comparison result under an additional
condition that the K$\ddot{a}$hler manifold $M^n$ is of maximum
volume growth. \vskip 0.5cm \noindent{\bf Theorem 1.1} ( Ni [14])
{\sl \ \ Let $M^n$ be a complete noncompact K$\ddot{a}$hler
manifold of complex dimension $n$ with nonnegative holomorphic
bisectional curvature. Assume that $M^n$ is of maximum volume
growth, i.e.,
$$Vol(B(x_0,r))\geq cr^{2n},\ \ \ for\ all\ 0\leq r<+\infty,$$
for some point $x_0\in M^n$ and some positive constant $c$. Then
$$dim_{\mathbb{C}}{\mathcal{O}}_d(M^n)\leq
dim_{\mathbb{C}}{\mathcal{O}}_{[d]}({\mathbb{C}}^n),\ \ \ for\
each\ d>0,$$ where $[d]$ is the greatest integer less than or
equal to $d$.}

Our first result in this paper is the following theorem which
resolves the remainder cases for the dimension comparison and
establishes the rigidity. \vskip 0.5cm\noindent {\bf Theorem 1.2}
{\sl \ \ Let $M^n$ be a complete noncompact K$\ddot{a}$hler
manifold of complex dimension $n$ with nonnegative holomorphic
bisectional curvature. Then there holds
$$dim_{\mathbb{C}}{\mathcal{O}}_d(M^n)\leq
dim_{\mathbb{C}}{\mathcal{O}}_{[d]}({\mathbb{C}}^n),\ \ \ for\
each\ d>0.$$ Moreover the equality holds for some positive integer
$d$ if and only if $M^n$ is holomorphically isometric to the
complex Euclidean space $\mathbb{C}^n$ with the standard flat
metric.}

We next try to improve the above dimension estimate when the
K$\ddot{a}$hler manifolds are not isometric to the flat complex
Euclidean space $\mathbb{C}^n$. The second result of this paper is
the following.\vskip 0.5cm\noindent {\bf Theorem 1.3} {\sl \ \ Let
$M^n$ be a complete noncompact K$\ddot{a}$hler manifold of complex
dimension $n$ with nonnegative holomorphic bisectional curvature.
Assume that the Ricci curvature of $M^n$ is positive at least at
one point in $M^n$. Then there exists a positive constant
$\epsilon\in(0,1)$ such that
$$dim_{\mathbb{C}}{\mathcal{O}}_d(M^n)\leq
dim_{\mathbb{C}}{\mathcal{O}}_{[(1-\epsilon)d]}({\mathbb{C}}^n),$$
for any positive integer $d>0$.}

We remark that the positive constant $\epsilon$ is depending on
the manifold. In Section 3, we will present a family of complete
noncompact K$\ddot{a}$hler manifolds with positive holomorphic
bisectional curvature such that for each $\epsilon\in(0,1)$, there
exists an element of this family such that the equality in the
above Theorem 1.3 holds for each positive integer $d$. We also
remark that the assumption that the Ricci curvature is positive at
least at one point in $M^n$ is almost necessary. In fact, by
considering the product manifold $\mathbb{C}^p\times M_2^{n-p}$
with $1\leq p\leq n$, we have
$dim_{\mathbb{C}}\mathcal{O}_1(\mathbb{C}^p\times M_2^{n-p})\geq
p+1$ but
$dim_{\mathbb{C}}\mathcal{O}_{[1-\epsilon]}(\mathbb{C}^n)=1$ for
any $\epsilon>0$.

The third result of this paper is to estimate the dimension in
terms of the volume growth and to establish the corresponding
rigidity. \vskip 0.5cm\noindent {\bf Theorem 1.4} {\sl \ \ Let
$M^n$ be a complete noncompact K$\ddot{a}$hler manifold of complex
dimension $n$ with nonnegative holomorphic bisectional curvature.
Suppose there exist $x_0\in M$, and $k>0$, $C>0$ such that,
$$Vol(B(x_0,r))\leq C(1+r)^{2k},\ \ \ for\ all\ r\geq 0.$$
Then there holds
$$dim_{\mathbb{C}}{\mathcal{O}}_d(M^n)\leq
dim_{\mathbb{C}}{\mathcal{O}}_{[d]}({\mathbb{C}}^{[k]}),\ \ \ for\
each\ d>0.$$ Moreover the equality holds for some positive integer
$d$ if and only if $M^n$ is holomorphically isometric to
$\mathbb{C}^{[k]}\times M_2^{n-[k]}$ for some complete
K$\ddot{a}$hler manifold $M_2^{n-[k]}$ of complex dimension
$n-[k]$ with nonnegative holomorphic bisectional curvature which
carries no nontrivial holomorphic functions of polynomial growth.}

Observe that the assumption on the curvature and the Bishop volume
comparison theorem assert that
$$Vol(B(x_0,r))\leq\omega_{2n}r^{2n},\ \ \ for\ all\ r\geq0,$$
where $\omega_{2n}$ is the volume of the unit ball in
$\mathbb{R}^{2n}$. On the other hand, a theorem of Calabi and Yau
[18] asserts that $Vol(B(x_0,r))$ must grow at least linearly.
Hence the constant $k$ in Theorem 1.4 must exist and satisfy
$1\leq k\leq n$. In particular, Theorem 1.2 is a special case of
Theorem 1.4. The real counterpart of Theorem 1.4 for harmonic
functions of polynomial growth had been studied by Li-Tam [11],
[12], Kasue [7], Li[8] and Cheeger-Colding-Minicozzi [2] for the
special cases $d=1$ or real dimension $m=2$. However the general
case remains as open questions (see for example Question 9.6 and
Question 9.7 in [10]).

We organize the paper as follows: in Section 2 we prove Theorem
1.2; in Section 3 we prove a slightly generalized version of
Theorem 1.3; finally in Section 4 we give a proof of Theorem 1.4.
\vskip 1cm

\section{Sharp Dimension Estimate and Rigidity}\indent

Let $M^n$ be a complete noncompact K$\ddot{a}$hler manifold of
complex dimension $n$. Fix a point $x_0\in M^n$. We call a
holomorphic function $f$ of polynomial growth if there exists
$d\geq0$ and $C$ depending on $x_0$, $d$ and $f$ such that
$$ |f(x)|\leq C(1+r^d(x,x_0)), \ \ \ \forall\ x\in M^n,$$
where $r(x,x_0)$ is the distance between $x$ and $x_0$. For a
holomorphic function $f$ of polynomial growth on $M^n$, we define
the degree of $f$ by
$$deg(f)=\inf \left\{ d \left|\begin{array}{c}
|f(x)|\leq C(1+r^d(x,x_0)),\ \ \forall\ x\in M^n,for\ \\[1mm]
some\ d\geq0\ and\ C=C(d,x_0,f)
\end{array}\right\}\right.
$$
It is clear that the definition of $deg(f)$ is independent of the
choice of the point $x_0$ in $M^n$. For any $d\geq0$ we denote by
$\mathcal{O}_d(M^n)$ the complex linear space of all holomorphic
functions of polynomial growth of degree at most $d$ on the
K$\ddot{a}$hler manifold $M^n$.

For any $f\in \mathcal{O}_d(M^n)$, we consider the following heat
equation
$$\left\{\begin{array}{l}
(\frac{\partial}{\partial t}-\Delta)u(x,t)=0,\ \ \ x\in M^n,\
t>0,\\[3mm]
u(x,0)=log|f(x)|^2,\ \ \ x\in M^n.
\end{array}\right.\eqno(2.1)
$$
Assuming the K$\ddot{a}$hler manifold $M^n$ has nonnegative Ricci
curvature, it is shown by Ni (see Lemma 3.1 in [14]) that the
equation (2.1) admits a smooth solution $u(x,t)$ on
$M^n\times(0,+\infty)$ and the solution $u(x,t)$ is given by
$$u(x,t)=\int_M H(x,y,t)log|f(y)|^2dy,\ \ \ \forall\ x\in M,\
t>0,\eqno(2.2)$$ where $H(x,y,t)$ is the heat kernel of $M^n$.
\vskip 0.5cm\noindent {\bf Lemma 2.1} {\sl \ \ Let $M^n$ be a
complete noncompact K$\ddot{a}$hler manifold with nonnegative
Ricci curvature and let u(x,t) solve the heat equation (2.1). Then
$$ \limsup_{t\rightarrow+\infty}\frac{u(x,t)}{log\ t}\leq d.$$}
\vskip 0.1cm\noindent{\sc Proof.} \ We write the solution $u(x,t)$
in terms of the heat kernel as
\begin{displaymath}
\begin{split}
    u(x,t)&=\int_M H(x,y,t)log|f(y)|^2dy\\[3mm]
          &=\int_{\{r(x,y)\leq\sqrt{t}\}}H(x,y,t)log|f(y)|^2dy
           +\int_{\{r(x,y)>\sqrt{t}\}}H(x,y,t)log|f(y)|^2dy,
\end{split}
\end{displaymath}
where $r(x,y)$ is the distance function from the point $x\in M^n$.
Since $f\in\mathcal{O}_d(M^n)$, for arbitrary $\epsilon>0$, there
exists a positive constant $C=C(d,\epsilon,x,f)$ such that
$$|f(y)|\leq C(1+r^{d+\epsilon}(x,y)),\ \ \ \forall\ y\in M^n.$$
Then we have for $t>1$,
\begin{displaymath}
\begin{split}
    u(x,t)&\leq\int_{\{r(x,y)\leq\sqrt{t}\}}H(x,y,t)((d+\epsilon)logt+C_1)dy
           \\[3mm]
           &\;\;\;+\int_{\{r(x,y)>\sqrt{t}\}}H(x,y,t)((d+\epsilon)log
           r^2(x,y)+C_2)dy\\[3mm]
          &\leq
          (d+\epsilon)logt+C_3+(d+\epsilon)\int_{\{r(x,y)>\sqrt{t}\}}H(x,y,t)log(\frac{r^2(x,y)}{t})dy
\end{split}
\end{displaymath}
for some positive constant $C_1$, $C_2$, $C_3$ depending only on
$d$, $\epsilon$, $x$ and $f$. Here we have used the fact that
$\int_{M^n}H(x,y,t)dy=1$.

By applying the heat kernel estimate of Li-Yau [13] and the
standard volume comparison, there exist positive constant $C(n)$
and $\tilde{C}(n)$ depending only on the dimension such that
\begin{displaymath}
\begin{split}
         &\int_{\{r(x,y)>\sqrt{t}\}}H(x,y,t)log(\frac{r^2(x,y)}{t})dy\\[3mm]
    \leq &C(n)\int_{\{r(x,y)>\sqrt{t}\}}\frac{1}{Vol(B(x,\sqrt{t}))}exp(-\frac{r^2(x,y)}{5t})\cdot
          log(\frac{r^2(x,y)}{t})dy\\[3mm]
    \leq &C(n)\sum_{k=0}^{\infty}\int\limits_{\{2^{k+1}\sqrt{t}\geq r(x,y)\geq2^k\sqrt{t}\}}
          \frac{1}{Vol(B(x,\sqrt{t}))}exp(-\frac{r^2(x,y)}{5t})\cdot log(\frac{r^2(x,y)}{t})dy\\[3mm]
    \leq &C(n)\sum_{k=0}^{\infty}exp(-\frac{2^{2k}}{5})\cdot log2^{2(k+1)}\cdot(2^{k+1})^{2n}\\[3mm]
    \leq &\tilde{C}(n).
\end{split}
\end{displaymath}
Hence we have for $t>1$,
$$u(x,t)\leq (d+\epsilon)logt+C_3+(d+\epsilon)\tilde{C}(n)$$
which completes the proof of Lemma 2.1. \hskip 1cm \#

For a nontrivial holomorphic function $f\in\mathcal{O}_d(M^n)$,
the vanishing order of $f$ at a fixed point $x\in M^n$ is defined
by
$$ ord_x(f)=\max\{m|\ D^\alpha f(x)=0,\ \forall\ |\alpha|<m\}.$$
The following result was first obtained by Ni in [14] for the case
that the K$\ddot{a}$hler manifold $M^n$ has nonnegative
holomorphic bisectional curvature. \vskip 0.5cm\noindent {\bf
Lemma 2.2}{\sl \ \ Let $M^n$ be a complete noncompact
K$\ddot{a}$hler manifold with nonnegative Ricci curvature. And let
$u(x,t),\ x\in M$ and $t\in(0,+\infty)$, be the solution of heat
equation (2.1) with $f\in \mathcal{O}_d(M^n)$ to be a nontrivial
holomorphic function. Denote by $w(x,t)=\frac{\partial}{\partial
t}u(x,t)$.Then
$$\lim_{t\rightarrow0}tw(x,t)=ord_x(f).$$}
\vskip 0.1cm\noindent {\sc Proof.} \ It is well known (see for
example [6]) that the vanishing order of a holomorphic function
$f$ at a fixed point $x\in M$ agrees with the Lelong number of the
associated positive (1,1) current
$\frac{\sqrt{-1}}{2\pi}\partial\bar{\partial}log(|f|^2)$ at the
point $x$, i.e.
$$ord_x(f)=\frac{1}{4n}\lim_{r\rightarrow0}\frac{r^2}{Vol(B(x,r))}\int_{B(x,r)}\Delta
log|f|^2dv.$$ Thus for any $\epsilon>0$, there exists a positive
constant $\delta<1$ such that
$$4n\ ord_x(f)-\epsilon\leq\frac{r^2}{Vol(B(x,r))}\int_{B(x,r)}\Delta
log|f|^2dv\leq 4n\ ord_x(f)+\epsilon,\eqno(2.3)$$ for all
$r\leq\delta$. Since $u(x,t)$ solves the heat equation (2.1) and
$w(x,t)=\frac{\partial}{\partial t}u(x,t)$, we see that for $t>0$,
\begin{displaymath}
\begin{split}
    tw(x,t)&=t\int_{M^n} H(x,y,t)\Delta log|f(y)|^2dy\\[3mm]
           &=t\int_{\{r(x,y)\leq\delta\}}H(x,y,t)\Delta log|f(y)|^2dy
             +t\int_{\{r(x,y)>\delta\}}H(x,y,t)\Delta
             log|f(y)|^2dy\\[3mm]
           &=I+II.
\end{split}
\end{displaymath}

To estimate the second term $II$, we claim that
$$\frac{R^2}{Vol(B(x,R))}\int_{B(x,R)}\Delta log|f|^2dy\leq
C_1log(R+2)+C_2,\\
\eqno(2.4)$$ for all $R>0$, where $C_1$, $C_2$ are positive
constants depending only on $n$, $d$, $x$ and $f$ (but independent
of $R$).

Choose a point $\tilde{x}$ near $x$ such that $f(\tilde{x})\neq0$.
Let $G_R$ be the positive Green's function on $B(\tilde{x},R)$
with Dirichlet boundary data. Then
\begin{displaymath}
\begin{split}
           &\int_{B(\tilde{x},R)}G_R(\tilde{x},y)\Delta log|f(y)|^2dy\\[3mm]
          =&-log|f(\tilde{x})|^2
            -\int_{\partial B(\tilde{x},R)}\frac{\partial G_R}{\partial \nu}\cdot log|f(y)|^2\\[3mm]
       \leq&-log|f(\tilde{x})|^2+Clog(R+2)
\end{split}
\end{displaymath}
for some positive constant $C$ depends only on $d$, $x$ and $f$.
Here we have used the facts that $\frac{\partial G_R}{\partial
\nu}<0$ on $\partial B(\tilde{x},R)$ and $\int_{\partial
B(\tilde{x},R)}\frac{\partial G_R}{\partial \nu}=-1$. It is
well-known (see for example Lemma 1.1 in [17]) that
$$ G_R(\tilde{x},y)\geq
C(n)\int_{r(\tilde{x},y)}^R\frac{t}{Vol(B(\tilde{x},t))}dt$$ for
all $y\in B(\tilde{x},\frac{1}{5}R)$ for some positive constant
$C(n)$ depending only on the dimension $n$.\\
Hence
\begin{displaymath}
\begin{split}
           &-log|f(\tilde{x})|^2+Clog(R+2)\\[3mm]
       \geq&C(n)\int_{0}^{\frac{R}{5}}(\int_t^R\frac{s}{Vol(B(\tilde{x},s))}ds)
            (\int_{\partial B(\tilde{x},t)}\Delta
            log|f|^2)dt\\[3mm]
          =&C(n)[(\int^{R}_{\frac{R}{5}}\frac{t}{Vol(B(\tilde{x},t))}dt)
            (\int_{B(\tilde{x},\frac{R}{5})}\Delta log|f|^2)
            +\int_0^{\frac{R}{5}}(\frac{t}{Vol(B(\tilde{x},t))}\int_{B(\tilde{x},t)}\Delta log|f|^2)dt]\\[3mm]
       \geq&\tilde{C}(n)\frac{R^2}{Vol(B(\tilde{x},\frac{R}{5}))}
            \int_{B(\tilde{x},\frac{R}{5})}\Delta log|f|^2,
\end{split}
\end{displaymath}
for some positive constant $\tilde{C}(n)$ depending only on the
dimension. Here we have used the standard volume comparison and
the fact that $\Delta log|f|^2\geq0$. This gives the claimation
(2.4).

By Li-Yau's upper bound estimate for the heat kernel [13] and
(2.4), the term $II$ can be estimated as
\begin{displaymath}
\begin{split}
         II&=t\int_{\{r(x,y)>\delta\}}H(x,y,t)\Delta log|f(y)|^2dy\\[3mm]
           &\leq C(n)t[ \int_{\{1\geq r(x,y)>\delta\}}\frac{1}{Vol(B(x,\sqrt{t}))}
             exp(-\frac{r^2(x,y)}{5t})\Delta log|f(y)|^2dy\\[3mm]
           & \ \ \ \ \ \ \ \ \ \ \ \
             +\sum_{k=1}^{\infty}\int\limits_{\{2^k\geq r(x,y)\geq2^{k-1}\}}\frac{1}{Vol(B(x,\sqrt{t}))}
             exp(-\frac{r^2(x,y)}{5t})\Delta log|f(y)|^2dy]\\[3mm]
           &\leq C(n)t[exp(-\frac{\delta^2}{5t})\cdot\frac{Vol(B(x,1))}{Vol(B(x,\sqrt{t}))}
              +C\sum_{k=1}^{\infty}exp(-\frac{2^{2(k-1)}}{5t})\cdot\frac{Vol(B(x,2^k))}{Vol(B(x,\sqrt{t}))}
               \cdot\frac{log(2^k+2)}{(2^k)^2}]\\[3mm]
           &\rightarrow0,\ \ \ as\ t\rightarrow0,
\end{split}
\end{displaymath}
$$\eqno(2.5)$$
for some positive constant $C$ depending only on $n$, $d$, $x$ and
$f$.

The following estimates for the term $I$,
$$\limsup_{t\rightarrow0}I\leq ord_x(f)+\epsilon,\eqno(2.6)$$
and$$ord_x(f)-\epsilon\leq\liminf_{t\rightarrow0}I.\eqno(2.7)$$
were obtained by Ni in [14] (in the proof of Lemma 4.1 in [14])
without using assumption on bisectional curvature. For the sake of
completeness, we present the proof as follows. It is well known
that the heat kernel
$$H(x,y,t)\sim\frac{1}{(4\pi t)^n}exp(-\frac{r^2(x,y)}{4t})+lower\
order\ term$$ as $t\rightarrow0$. Then for $t>0$ sufficiently
small,
\begin{displaymath}
\begin{split}
         I&\leq t\int_{0}^{\delta}\frac{1}{(4\pi t)^n}exp(-\frac{s^2}{4t})
           (\int_{\partial B(x,s)}\Delta log|f(y)|^2dy)ds+\frac{\epsilon}{2}\\[3mm]
          &=\frac{t}{(4\pi t)^n}exp(-\frac{\delta^2}{4t})
           (\int_{B(x,\delta)}\Delta log|f(y)|^2dy)\\[3mm]
          & \ \ \ +t\int_{0}^{\delta}\frac{1}{(4\pi t)^n}exp(-\frac{s^2}{4t})
           (\int_{B(x,s)}\Delta log|f(y)|^2dy)(\frac{s}{2t})ds+\frac{\epsilon}{2}\\[3mm]
          &=III+IV+\frac{\epsilon}{2}
\end{split}
\end{displaymath}
Clearly the term $III$ has limit $0$ as $t\rightarrow0$. The term
$IV$ can be estimated as
\begin{displaymath}
\begin{split}
        IV&\leq t\int_{0}^{\delta}\frac{\omega_{2n}s^{2n}}{(4\pi t)^n}exp(-\frac{s^2}{4t})
           (\frac{1}{Vol(B(x,s))}\int_{B(x,s)}\Delta log|f(y)|^2dy)(\frac{s}{2t})ds\\[3mm]
          &=\frac{\omega_{2n}}{\pi^n}\int_{0}^{\frac{\delta^2}{4t}}\frac{1}{4}exp(-\tau)\tau^{n-1}
           (\frac{4t\tau}{Vol(B(x,\sqrt{4t\tau}))}\int_{B(x,\sqrt{4t\tau})}\Delta log|f(y)|^2dy)d\tau\\[3mm]
          &\leq n\frac{\omega_{2n}}{\pi^n}ord_x(f)\int_{0}^{\frac{\delta^2}{4t}}
           exp(-\tau)\tau^{n-1}d\tau+\frac{\epsilon}{2}\\[3mm]
          &=ord_x(f)+\frac{\epsilon}{2}
\end{split}
\end{displaymath}
by using the standard volume comparison and (2.3). This proves
(2.6). The proof for (2.7) is similar.

Therefore the combination of (2.5), (2.6) and (2.7) completes the
proof of Lemma 2.2.\# \vskip 0.5cm We are now ready to prove
Theorem 1.2.

\noindent{\bf Proof of Theorem 1.2}

Let $M^n$ be a complete noncompact K$\ddot{a}$hler manifold of
complex dimension $n$ with nonnegative holomorphic bisectional
curvature. For any nontrivial $f\in\mathcal{O}_d(M^n)$, let
$u(x,t)$, on $M^n\times(0,+\infty)$, be the solution of the heat
equation (2.1). When we consider the trivial line bundle
$M^n\times\mathbb{C}\rightarrow M^n$ equipped with the metric
$e^{-u(x,t)}$ on the fibers, the complex Hessian $u_{i\bar{j}}$
corresponds the bundle curvature tensor and the heat equation
(2.1) is exactly the Hermitian Yang-Mills flow for Hermitian
metric.

Differentiate the equation (2.1) we have for normal coordinates at
a point
$$(u_t)_{i\bar{j}}=R_{l\bar{k}i\bar{j}}u_{k\bar{l}}+\frac{\partial^4u}{\partial
z^k\partial \bar{z}^k\partial z^i\partial \bar{z}^j}$$ and by
definition and an easy computation
\begin{displaymath}
\begin{split}
    \Delta u_{i\bar{j}}&=\frac{1}{2}(u_{i\bar{j},k\bar{k}}+u_{i\bar{j},\bar{k}k})\\[3mm]
                      &=\frac{\partial^4u}{\partial z^i\partial \bar{z}^j\partial z^k\partial \bar{z}^k}
                       +\frac{1}{2}(R_{i\bar{l}}u_{l\bar{j}}+R_{l\bar{j}}u_{i\bar{l}}).
\end{split}
\end{displaymath}
Thus the bundle curvature tensor $u_{i\bar{j}}$ satisfies the
complex Lichnerowicz-Laplacian heat equation
$$(\frac{\partial}{\partial t}-\Delta)u_{i\bar{j}}=R_{l\bar{k}i\bar{j}}u_{k\bar{l}}
-\frac{1}{2}(R_{i\bar{l}}u_{l\bar{j}}+R_{l\bar{j}}u_{i\bar{l}}).\eqno(2.8)$$
It was shown by Ni [14] (Lemma 2.1 in [14]) that the nonnegativity
of the bundle curvature tensor $u_{i\bar{j}}$ is preserved for all
$t>0$. More importantly, Ni [14] (Corollary 2.2 in [14]) proved
that the Hermitian Yang-Mills flow possesses a Li-Yau-Hamilton
inequality (as Ricci flow or mean curvature flow)
$$Z(V)\stackrel{\Delta}{=}w_t+\nabla_kw\cdot V^k+\nabla_{\bar{k}}w\cdot
V^{\bar{k}}+u_{i\bar{j}}V^iV^{\bar{j}}+\frac{w}{t}\geq0,\eqno(2.9)$$
on $M^n\times(0,+\infty)$ for any (1,0) vector field $V$. Here
$w(x,t)=\frac{\partial}{\partial t}u(x,t)=\Delta u(x,t)$ is the
trace of the bundle curvature tensor. In particular by choosing
$V\equiv0$,
$$\frac{\partial}{\partial t}(tw(x,t))\geq0,\eqno (2.10)$$
on $M^n\times(0,+\infty)$. Then for a large enough positive
integer $p$, we have
\begin{displaymath}
\begin{split}
      u(x,t)=\int_1^tw(x,s)ds+u(x,1)&\geq\int_{t^{\frac{1}{p}}}^t\frac{t^{\frac{1}{p}}}{s}
                                     w(x,t^{\frac{1}{p}})ds+u(x,1)\\[3mm]
                                    &=(1-\frac{1}{p})logt\cdot
                                     t^{\frac{1}{p}}w(x,t^{\frac{1}{p}})+u(x,1).
\end{split}
\end{displaymath}
By combining with Lemma 2.1, we have
$$\limsup_{t\rightarrow\infty}(1-\frac{1}{p})t^\frac{1}{p}w(x,t^\frac{1}{p})\leq
d$$ and by arbitrariness of $p$,
$$\limsup_{t\rightarrow\infty}tw(x,t)\leq d$$
So by combining with (2.10), we deduce
$$tw(x,t)\leq d,\ \ \ for\ all\ t\geq0.\eqno(2.11)$$
On the other hand, from Lemma 2.2,
$$\lim_{t\rightarrow0}tw(x,t)=ord_x(f).$$
Hence
$$ ord_x(f)\leq[d],\ \ \ for\ any\ nontrivial\
f\in\mathcal{O}_d(M^n),\eqno(2.12)$$ since the vanishing order
need to be integer.

In the following we perform the standard Poincar$\acute{e}$-Siegel
argument. For the fixed $x\in M^n$ and a local holomorphic
coordinate system $\{z^1,\cdots,z^n\}$ near $x$, define the
Poincar$\acute{e}$-Siegel map
$$
\arraycolsep=1.5pt
\begin{array}{lccl}
P:&\mathcal{O}_d(M^n)&\rightarrow&\mathbb{C}^{k_{[d]}}\\[3mm]
&f&\longmapsto&(f(x),Df(x),\cdots,D^\alpha f(x),\cdots),\ \ \ for\
all\ |\alpha|\leq[d]
\end{array}
$$
where $k_{[d]}=dim_{\mathbb{C}}(\mathcal{O}_{[d]}(\mathbb{C}^n))$.
The inequality (2.12) implies that the map $P$ must be injective.
Thus we have proved
$$dim_{\mathbb{C}}\mathcal{O}_d(M^n)\leq
dim_{\mathbb{C}}\mathcal{O}_{[d]}(\mathbb{C}^n)\eqno(2.13)$$

We now study the rigidity part of Theorem 1.2. Suppose there
exists a positive integer $d$ such that
$$dim_{\mathbb{C}}\mathcal{O}_d(M^n)=dim_{\mathbb{C}}\mathcal{O}_{d}(\mathbb{C}^n)$$
From the above argument we know that there exists a nonconstant
function $f\in\mathcal{O}_d(M^n)$ and a point $x_0\in M^n$ such
that
$$ord_{x_0}(f)=d$$
Let $u(x,t)$ be the solution of the heat equation (2.1) with the
initial data $log|f|^2$, and $w(x,t)=\frac{\partial}{\partial
t}u(x,t)=\Delta u(x,t)$. By (2.11) and Lemma 2.2, we have
$$tw(x_0,t)\equiv d,\ \ \ for\ all\ t>0,$$
which implies
$$Z(0)=(w_t+\frac{w}{t})(x_0,t)\equiv0,\ \ \ for\ all\
t>0.\eqno(2.14)$$

In the following we distinguish two cases. Case(i):
$(u_{i\bar{j}}(x,t))$ is positive definite on
$M^n\times(0,+\infty)$ and Case(ii): $(u_{i\bar{j}}(x,t))$ has
nontrivial kernel for some $(x,t)\in M^n\times(0,+\infty)$.

Let us first consider Case(i). Recall from (2.9) that
$$Z(V)\stackrel{\Delta}{=}w_t+\frac{w}{t}+\nabla_kw\cdot V^k+\nabla_{\bar{k}}w\cdot
V^{\bar{k}}+u_{i\bar{j}}V^iV^{\bar{j}}\geq0$$ on
$M^n\times(0,+\infty)$ for any (1,0) vector field $V$. Since
$(u_{i\bar{j}})>0$, we see that there exists a unique smooth
vector field $V$ which minimizes $Z$ and satisfies
$$\nabla_kw+u_{k\bar{j}}V^{\bar{j}}=0,\ \ \ everywhere,$$
by the first variation formula. We then consider the function
$Z=Z(V)$ by taking $V$ to be the unique minimizing vector field.
It was shown by Ni in Lemma 1.3 of [14] that the function $Z$
satisfies, for normal coordinates at a point,

\begin{displaymath}
\begin{split}
    (\frac{\partial}{\partial t}-\Delta)Z=&R_{i\bar{j}k\bar{l}}u_{l\bar{k}}V^iV^{\bar{j}}+u_{i\bar{j}}
               (\nabla_{\bar{k}}V_j-\frac{1}{t}g_{j\bar{k}})(\nabla_{k}V_{\bar{i}}-\frac{1}{t}g_{k\bar{i}})\\[3mm]
               &+u_{i\bar{j}}\nabla_{\bar{k}}V_{\bar{i}}\nabla_kV_j-\frac{2Z}{t}.
\end{split}
\end{displaymath}
$$\eqno(2.15)$$
Note that the bisectional curvature is nonnegative, $u$ is
plurisubharmonic and by (2.14) $Z(x_0,t)=0$ for all $t>0$. Thus by
the strong maximum principle, we see that the function $Z$ must be
identically zero everywhere. And from (2.15) we have
$$ R_{i\bar{j}k\bar{l}}u_{l\bar{k}}V^iV^{\bar{j}}=0,\eqno(2.16)$$
$$ u_{i\bar{j}}(\nabla_{\bar{k}}V_j-\frac{1}{t}g_{j\bar{k}})
(\nabla_{k}V_{\bar{i}}-\frac{1}{t}g_{k\bar{i}})=0,\eqno(2.17)$$
$$u_{i\bar{j}}\nabla_{\bar{k}}V_{\bar{i}}\nabla_kV_j=0,\eqno(2.18)$$
everywhere on $M^n\times(0,+\infty)$. Since we assumed that
$(u_{i\bar{j}}(x,t))$ is positive definite everywhere in Case(i),
we have from (2.17) and (2.18),
$$\left\{
\begin{array}{l}
  \nabla_iV_{\bar{j}}=\nabla_{\bar{j}}V_i=\frac{1}{t}g_{i\bar{j}}\\[3mm]
  \nabla_iV_j=\nabla_{\bar{i}}V_{\bar{j}}=0.
\end{array}\right.
$$
Let us first assume that the manifold $M^n$ is simply connected.
Then there exists a function $\varphi$ on $M^n$ such that
$$\left\{
\begin{array}{l}
  \varphi_{i\bar{j}}=g_{i\bar{j}}\\[3mm]
  \nabla_i\nabla_j\varphi=\nabla_{\bar{i}}\nabla_{\bar{j}}\varphi=0.
\end{array}\right.\eqno(2.19)
$$
Note that for any normal geodesic $\gamma$, we have
$$\frac{d^2}{ds^2}\varphi\circ\gamma(s)=Hess(\varphi)(\gamma',\gamma')=1.\eqno(2.20)$$
This implies that $\varphi$ is a proper convex function. Let
$p_0\in M^n$ be the unique minimizer of $\varphi$. We may assume
$\varphi(p_0)=0$. Then from (2.20) it is clear that
$$\varphi(x)=\frac{1}{2}r^2(x,p_0),\ \ \ on\ M^n.$$
Hence by (2.19) and the rigidity part of standard Hessian
comparison theorem (or Laplacian comparison theorem), we know that
$M^n$ is flat and then $M$ is holomorphically isometric to
$\mathbb{C}^n$ since $M^n$ is assumed to be simply connected.

To remove the simply connectedness assumption, we note that
$$dim_{\mathbb{C}}\mathcal{O}_d(M^n)\leq
dim_{\mathbb{C}}\mathcal{O}_d(\tilde{M}^n),$$ where $\tilde{M}^n$
is the universal cover of $M^n$. Thus from (2.13) and the
assumption on $d$, we have
$$dim_{\mathbb{C}}\mathcal{O}_d(\tilde{M}^n)=
dim_{\mathbb{C}}\mathcal{O}_d(\mathbb{C}^n),$$ and then by the
above argument, $\tilde{M}^n$ is holomorphically isometric to
$\mathbb{C}^n$. Note that every function in $\mathcal{O}_d(M^n)$
can be lifted as a function of $\mathcal{O}_d(\tilde{M}^n)$. This
says, the complex linear space $\mathcal{O}_d(M^n)$ is a subspace
of $\mathcal{O}_d(\tilde{M}^n)$. Since
$$dim_{\mathbb{C}}\mathcal{O}_d(M^n)=dim_{\mathbb{C}}\mathcal{O}_d(\tilde{M}^n)=
dim_{\mathbb{C}}\mathcal{O}_d(\mathbb{C}^n),$$ we see that every
function in $\mathcal{O}_d(\tilde{M}^n)$ is
$\pi_1(M^n)$-invariant. Let $z^1,\cdots,z^n$ be the coordinate
functions on $\tilde{M}^n=\mathbb{C}^n$. Obviously
$z^1,\cdots,z^n\in\mathcal{O}_d(\tilde{M}^n)$. Let
$\sigma\in\pi_1(M^n)$ be a deck transformation. Since the
functions $z^1,\cdots,z^n$ are $\pi_1(M^n)$-invariant, we have
$$\sigma([z^i=0])=[z^i=0],\ \ \ for\ i=1,\cdots,n.$$
Hence
$$\sigma(\{0\})=\sigma(\bigcap\limits_{i=1}^n[z^i=0])=\bigcap\limits_{i=1}^n[z^i=0]=\{0\}.$$
This implies $\pi_1(M^n)=0$, so we finish the proof of Case(i).

We next consider Case(ii). The argument of the previous paragraph
tells us that we may assume $M^n$ is simply connected. Recall that
$(u_{i\bar{j}}(x,t))$ is nonnegative definite on
$M^n\times(0,+\infty)$ and satisfies the complex
Lichnerowicz-Laplacian heat equation (2.8). Suppose that the
matrix $(u_{i\bar{j}}(x,t))$ has nontrivial kernel at some
$(\bar{x},\bar{t})\in M^n\times (0,+\infty)$. Then by (2.8) and
the strong maximum principle, we know that for all $t<\bar{t}$ and
$x\in M^n$, $(u_{i\bar{j}}(x,t))$ has nontrivial kernel. Denote
the kernel space of $(u_{i\bar{j}}(x,t))$ by $K(x,t)\subset
T_x^{1,0}M^n$. It was shown by Ni-Tam in Corollary 2.1 of [15]
that there exists $\tilde{t}\in(0,\bar{t})$ such that for any
$0<t<\tilde{t}$, $K(x,t)$ is a distribution which is invariant
under parallel translations. Moreover the K$\ddot{a}$hler manifold
$M^n$ splits isometrically and holomorphically as
$$M^n=M^p_1\times M^{n-p}_2$$
with $1\leq p\leq n$, where $K$ corresponds the tangent bundle of
$M_1^p$ and $(u_{i\bar{j}}(x,t))>0$ on
$M_2^{n-p}\times(0,\tilde{t})$. Both $M_1^p$ and $M_2^{n-p}$ are
complete K$\ddot{a}$hler manifold with nonnegative holomorphic
bisectional curvature.

Write $x=(x_1,x_2)\in M^n=M_1^p\times M_2^{n-p}$ with $x_1\in
M_1^p$ and $x_2\in M_2^{n-p}$. Since the restriction
$u_{i\bar{j}}|_{M_1^p\times(0,\tilde{t})}\equiv0$, we have that
for any fixed $x_2\in M_2^{n-p}$ and $t\in(0,\tilde{t})$, the
function $u(\cdot,x_2,t)$ is a pluriharmonic function on $M_1^p$.
This is, $$d_{x_1}d_{x_1}^cu(\cdot,x_2,t)=0,\ \ \ on\ M_1^p,$$
where $d_{x_1}$ is the exterior differential on $M_1^p$ and
$d_{x_1}^c=\sqrt{-1}(\bar{\partial}_{x_1}-{\partial}_{x_1})$ is
the usual real operator. By the heat kernel estimate of Li-Yau
[13], it is not hard (see for example Corollary 1.4 in [15]) to
see that the function $u(\cdot,x_2,t)$ is at most upper
logarithmic growth on $M_1^p$ (but we can not insure such a lower
bound). Since we assumed $M^n$ is simply connected, $M_1^p$ is
also simply connected. Thus there exists a real function $v$ on
$M_1^p$ such that $$d^c_{x_1}u(\cdot,x_2,t)=d_{x_1}v(\cdot),\ \ \
on\ M_1^p,$$ i.e.,
$$\sqrt{-1}(\bar{\partial}_{x_1}-{\partial}_{x_1})u(\cdot,x_2,t)=(\partial_{x_1}+\bar{\partial}_{x_1})v(\cdot),\
\ \ on\ M_1^p,$$ which implies
$$\bar{\partial}_{x_1}(u(\cdot,x_2,t)+\sqrt{-1}v(\cdot))=0,\ \ \ on\ M_1^p.$$
So $u(\cdot,x_2,t)+\sqrt{-1}v(\cdot)$ is a holomorphic function on
$M_1^p$. Let $\alpha>0$ be a positive constant and set
$$F(\cdot)=e^{\alpha(u(\cdot,x_2,t)+\sqrt{-1}v(\cdot))},\ \ \ on\ M_1^p.$$
Clearly the function $F$ is holomorphic on $M_1^p$. When we choose
$\alpha>0$ small enough, there holds $$|F(x_1)|\leq
C(1+r(x_1,x_0^{(1)}))^{\frac{1}{2}},$$ for all $x_1\in M_1^p$ and
for some fixed $x_0^{(1)}\in M_1^p$, since the function
$u(\cdot,x_2,t)$ is at most upper logarithmic growth on $M_1^p$.
It then follows from Cheng-Yau [4] (or the vanishing order
estimate (2.12)) that $F$ is constant on $M_1^p$. This implies
that $u(\cdot,x_2,t)$ is constant on $M_1^p$. Hence the solution
$u(x,t)$ can be regarded as a solution of the heat equation (2.1)
on $M_2^{n-p}\times(0,\tilde{t})$ with positive definite
$(u_{i\bar{j}}(x,t))$. Clearly the solution $u(x,t)$ is actually
defined on $M_2^{n-p}\times(0,+\infty)$ and $(u_{i\bar{j}}(x,t))$
is positive definite everywhere (by the strong maximum principle)
on $M_2^{n-p}\times(0,+\infty)$.

The equation (2.14) tells us the function
$w(x,t)=\frac{\partial}{\partial t}u(x,t)=\Delta u(x,t)$
satisfying
$$Z(0)=(w_t+\frac{w}{t})(x_0,t)\equiv0,\ \ \ for\ all\ t>0,$$
for some $x_0\in M_2^{n-p}$. Then by repeating the argument of
Case(i) we conclude that $M_2^{n-p}$ is holomorphically isometric
to $\mathbb{C}^{n-p}$. So
$$M^n=M_1^p\times\mathbb{C}^{n-p},\eqno(2.21)$$
isometrically and holomorphically. Since the function
$f\in\mathcal{O}_d(M^n)$ with $ord_{x_0}(f)=d>0$ is nonconstant,
it is clear from the strong maximum principle (see for example
Lemma 3.1 in [14]) that $w(x,t)>0$ everywhere on
$M^n\times(0,+\infty)$. This says the rank of
$(u_{i\bar{j}}(x,t))$ is at least $1$. Hence we must have $p\leq
n-1$.

We now claim that
$$dim_{\mathbb{C}}\mathcal{O}_d(M_1^p)=
dim_{\mathbb{C}}\mathcal{O}_d(\mathbb{C}^p).\eqno(2.22)$$ Clearly
once this is proved, then by induction on the dimension of the
manifolds we will complete the proof of the rigidity part of
Theorem 1.2.

Fix a point $(p_0,q_0)\in M_1^p\times\mathbb{C}^{n-p}$ and a local
holomorphic coordinate system $(z_1,z_2)$ near the point
$(p_0,q_0)$, where $z_1=(z^1_1,\cdots,z^p_1)\in M_1^p$ and
$z_2=(z^1_2,\cdots,z^{n-p}_2)\in \mathbb{C}^{n-p}$. For any
holomorphic function $g\in\mathcal{O}_d(M^n)$, $g$ has a Taylor
expansion near $(p_0,q_0)$. Let us denote $P_d(g)(z_1,z_2)$ be the
polynomial obtained by truncating the Taylor expansion up to order
$d$, i.e.,
$$g(z_1,z_2)=P_d(g)(z_1,z_2)+higher\ order\ terms.$$
Consider the Poincar$\acute{e}$-Siegel map
$$\arraycolsep=1.5pt\begin{array}{lccl}
P_d:&\mathcal{O}_d(M^n)&\rightarrow&\mathcal{O}_d(\mathbb{C}^n)
    =\mathcal{O}_d(\mathbb{C}^p\times\mathbb{C}^{n-p}),\\[3mm]
&g&\longmapsto&P_d(g).
\end{array}$$
The estimate (2.12) implies that the map $P_d$ is injective. Since
$$dim_{\mathbb{C}}\mathcal{O}_d(M^n)=dim_{\mathbb{C}}\mathcal{O}_d(\mathbb{C}^n)$$
by the assumption, we see that the Poincar$\acute{e}$-Siegel map
$P_d$ is an isomorphism between $\mathcal{O}_d(M^n)$ and
$\mathcal{O}_d(\mathbb{C}^n)$.

Since $M_1^p$ is also complete K$\ddot{a}$hler manifold with
nonnegative holomorphic bisectional curvature, we then also have
the estimate (2.12) for any nontrivial $f\in\mathcal{O}_d(M_1^p)$.
Thus the corresponding Poincar$\acute{e}$-Siegel map (by
considering the Taylor expansion at $p_0\in M_1^p$)
$$P_d^{(1)}:\mathcal{O}_d(M_1^p)\rightarrow
\mathcal{O}_d(\mathbb{C}^p)$$ is also injective. To prove the
claim (2.22), it suffices to prove $P_d^{(1)}$ is surjective. Of
course, we can regard $\mathcal{O}_d(M_1^p)$ as a complex subspace
of $\mathcal{O}_d(M^n)$ and this Poincar$\acute{e}$-Siegel map
$P_d^{(1)}$ is just the restriction of original one $P_d$ on the
subspace $\mathcal{O}_d(M_1^p)(\subset\mathcal{O}_d(M^n))$. Since
$P_d$ is a bijection between $\mathcal{O}_d(M^n)$ and
$\mathcal{O}_d(\mathbb{C}^n)$, we only need to prove the assertion
that if $g\in\mathcal{O}_d(M^n)$ with
$P_d(g)\in\mathcal{O}_d(\mathbb{C}^p)(\subset\mathcal{O}_d(\mathbb{C}^n))$,
then $g\in\mathcal{O}_d(M_1^p)$. In fact, for such a $g$ with
$P_d(g)\in\mathcal{O}_d(\mathbb{C}^p)$, $g$ has a Taylor expansion
of the form
$$g(z_1,z_2)=P_d(g)(z_1)+higher\ order\ terms$$
which implies that
$$D_{z_2}^{\alpha}g(z_1,z_2)|_{(p_0,q_0)}=0,\ \ \ for\ any\
|\alpha|\leq d\ and\ |\alpha|\neq0.\eqno(2.23)$$ Fix $z_1\in
M_1^p$ and consider the function $g(z_1,\cdot)$ on
$\mathbb{C}^{n-p}$. Clearly the function $g(z_1,\cdot)$ belongs to
$\mathcal{O}_d(\mathbb{C}^{n-p})$ since $g$ belongs to
$\mathcal{O}_d(M^n)$. Thus there also holds the vanishing order
estimate (2.12). So by combining (2.23) we see that $g(z_1,\cdot)$
is constant on $\mathbb{C}^{n-p}$. This says $g$ is only a
function of variable $z_1\in M_1^p$ and then
$g\in\mathcal{O}_d(M_1^p)$.

Therefore we have proved the claimation (2.22) and then have
completed the proof of Theorem 1.2.\hskip 1cm \# \vskip 1cm

\section{Improved Dimension Estimate}\indent

Let $M^n$ be a complete noncompact K$\ddot{a}$hler manifold of
complex dimension $n$ with nonnegative holomorphic bisectional
curvature. By inspecting the proof of the rigidity part of Theorem
1.2, we have actually proved that if there exists a holomorphic
function $f$ of polynomial growth on $M^n$ such that
$$ord_x(f)=deg(f)>0,\ \ \ for\ some\ x\in M^n,$$
then the universal cover $\tilde{M}^n$ admits a splitting
$$\tilde{M}^n=M_1^p\times\mathbb{C}^{n-p}$$
isometrically and holomorphically for some $p\leq n-1$. Thus if
$M^n$ does not admit its universal cover $\tilde{M}^n$ splitting
isometrically and holomorphically as $\mathbb{C}\times M_2^{n-1}$,
then we have
$$ord_x(f)\leq deg(f)-1,$$
for all $f\in \mathcal{O}_d(M^n)$, all points $x\in M^n$ and all
positive integers $d$, consequently,
$$dim_{\mathbb{C}}\mathcal{O}_d(M^n)\leq
dim_{\mathbb{C}}\mathcal{O}_{d-1}(\mathbb{C}^n),\ \ \ for\ all\
positive\ integers\ d.\eqno(3.1)$$ The following result improves
this dimension estimate significantly as $d$ large. \vskip
0.5cm\noindent {\bf Theorem 3.1}{\sl \ \ Let $M^n$ be a complete
noncompact K$\ddot{a}$hler manifold of complex dimension $n$ with
nonnegative holomorphic bisectional curvature. Suppose its
universal cover $\tilde{M}^n$ does not admit a holomorphically
isometric splitting as $\mathbb{C}\times M_2^{n-1}$. Then there
exists a positive constant $\epsilon\in(0,1)$ such that
$$dim_{\mathbb{C}}\mathcal{O}_d(M^n)\leq
dim_{\mathbb{C}}\mathcal{O}_{[(1-\epsilon)d]}(\mathbb{C}^n),\ \ \
for\ all\ positive\ integers\ d.$$} \vskip 0.1cm Before giving the
proof of Theorem 3.1, we remark that Theorem 1.3 in the
introduction is a consequence of the above Theorem 3.1. Indeed if
a K$\ddot{a}$hler manifold $M^n$ admits its universal cover
$\tilde{M}^n$ isometrically splitting as $\mathbb{C}\times
M_2^{n-1}$, then its Ricci curvature must have a nontrivial kernel
everywhere. \vskip 0.1cm \noindent{\bf Proof of Theorem 3.1}

We argue by contradiction. Suppose for each positive integer $k$,
there exists a positive integer $d_k$ such that
$$dim_{\mathbb{C}}\mathcal{O}_{d_k}(M^n)>
dim_{\mathbb{C}}\mathcal{O}_{[(1-\frac{1}{k})d_k]}(\mathbb{C}^n).$$
Since the universal cover $\tilde{M}^n$ does not split
isometrically and holomorphically as $\mathbb{C}\times M_2^{n-1}$,
we have the estimate (3.1). It follows that $d_k\geq k$ and then
$$\lim_{k\rightarrow \infty}\frac{[(1-\frac{1}{k})d_k]}{d_k}=1.$$
Fix a point $x_0\in M^n$. By the Poincar$\acute{e}$-Siegel
argument, we know that for each positive integer $k$ there exists
a function $f_k\in\mathcal{O}_{d_k}(M^n)$ such that
$$[(1-\frac{1}{k})d_k]\leq ord_{x_0}(f_k).$$
On the other hand, by (2.12) we always have
$$ord_{x_0}(f_k)\leq deg(f_k)=d_k,\ \ \ for\ each\ positive\
integer\ k.$$ Thus
$$\lim_{k\rightarrow\infty}\frac{ord_{x_0}(f_k)}{deg(f_k)}=1.$$

Let $u_k(x,t)$ solves the following heat equation
$$\left\{
\begin{array}{l}
 \frac{\partial}{\partial t}u_k=\Delta u_k,\ \ \ on\
  M^n\times(0,+\infty),\\[3mm]
 u_k|_{t=0}=log|f_k|^2,\ \ \ on\ M^n,
\end{array}\right.
$$
and let $w_k(x,t)=\Delta u_k(x,t)$. We have seen in Lemma 2.2 and
(2.11), for each $k$,
$$tw_k(x_0,t)\geq[(1-\frac{1}{k})deg(f_k)],$$
and
$$tw_k(x,t)\leq deg(f_k),$$
for all $t>0$ and $x\in M^n$.

Set
$$(\tilde{u}_k)_{i\bar{j}}=\frac{1}{deg(f_k)}(u_k)_{i\bar{j}}\ \ \
\ \ \ \ and\ \ \ \ \ \ \ \tilde{w}_k=\frac{1}{deg(f_k)}w_k.$$ Then
we have from (2.8),
$$\left\{\begin{array}{l}
\frac{\partial}{\partial
t}(\tilde{u}_k)_{i\bar{j}}=\Delta(\tilde{u}_k)_{i\bar{j}}+R_{l\bar{m}i\bar{j}}(\tilde{u}_k)_{m\bar{l}}
-\frac{1}{2}(R_{i\bar{l}}(\tilde{u}_k)_{l\bar{j}}
+R_{l\bar{j}}(\tilde{u}_k)_{i\bar{l}}),\ \ \ on\
M^n\times(0,+\infty),\\[4mm]
t\tilde{w}_k(x,t)\leq 1,\ \ \ on\ M^n\times(0,+\infty),\\[4mm]
t\tilde{w}_k(x_0,t)\geq\frac{[(1-\frac{1}{k})d_k]}{d_k}\rightarrow
1,\ \ \ as\ k\rightarrow+\infty,\ for\ all\ t>0,\\[4mm]
g^{i\bar{j}}(\tilde{u}_k)_{i\bar{j}}=\tilde{w}_k\ \ \ and\ \ \
((\tilde{u}_k)_{i\bar{j}})\ is\ nonnegative\ definite\ on\
M^n\times(0,+\infty).
\end{array}\right.$$
Since we have uniform bound of $(\tilde{u}_k)_{i\bar{j}}$ on every
compact subset of $M^n\times(0,+\infty)$, the standard Bernstein
trick for linear parabolic equations gives all the higher
derivative uniform estimates of $(\tilde{u}_k)_{i\bar{j}}$ on
every compact subset of $M^n\times(0,+\infty)$. Then there exists
a subsequence of $k\rightarrow\infty$ so that
$$(\tilde{u}_k)_{i\bar{j}}\rightarrow \tilde{u}_{i\bar{j}}$$
$$\tilde{w}_k\rightarrow\tilde{w}$$
for some smooth nonnegative (1,1)-form $\tilde{u}_{i\bar{j}}$ and
some smooth function $\tilde{w}$ in $C^\infty$ topology over
compact sets of $M^n\times(0,+\infty)$, which satisfy
$$\left\{\begin{array}{l}
\frac{\partial}{\partial
t}\tilde{u}_{i\bar{j}}=\Delta\tilde{u}_{i\bar{j}}+R_{l\bar{k}i\bar{j}}\tilde{u}_{k\bar{l}}
-\frac{1}{2}(R_{i\bar{l}}\tilde{u}_{l\bar{j}}
+R_{l\bar{j}}\tilde{u}_{i\bar{l}}),\ \ \ on\
M^n\times(0,+\infty),\\[4mm]
t\tilde{w}(x,t)\leq 1,\ \ \ on\ M^n\times(0,+\infty),\\[4mm]
t\tilde{w}(x_0,t)\geq1,\ \ \ for\ all\ t>0,\\[4mm]
g^{i\bar{j}}\tilde{u}_{i\bar{j}}=\tilde{w}\ \ \ and\ \ \
(\tilde{u}_{i\bar{j}})\ is\ nonnegative\ definite\ on\
M^n\times(0,+\infty).
\end{array}\right.$$
In particular, we have
$$t\tilde{w}(x_0,t)\equiv1,\ \ \ for\ all\ t>0.$$
Clearly we still has the Li-Yau-Hamilton inequality (2.9) as
$$Z(V)=\tilde{w}_t+\nabla_k\tilde{w}\cdot V^k+\nabla_{\bar{k}}\tilde{w}\cdot
V^{\bar{k}}+\tilde{u}_{i\bar{j}}V^iV^{\bar{j}}+\frac{\tilde{w}}{t}\geq0$$
on $M^n\times(0,+\infty)$ for any (1,0) vector field $V$. Thus the
same argument in the proof of the rigidity part of Theorem 1.2
works in the present situation. Therefore we conclude that the
universal cover $\tilde{M}^n$ splits as
$$\tilde{M}^n=M_1^p\times\mathbb{C}^{n-p}$$
isometrically and holomorphically for some $p\leq n-1$, which
contradicts with the assumption.\hskip 1cm\# \vskip 0.5cm A
natural question is to ask if we can choose the positive constant
$\epsilon\in(0,1)$ in Theorem 3.1 independent of the
K$\ddot{a}$hler manifold. We now present a family of examples to
show that it is impossible. In [1], Cao constructed a family of
complete K$\ddot{a}$hler metrics of the form
$$g^\lambda_{i\bar{j}}=\partial_i\partial_{\bar{j}}u_\lambda(log|z|^2),\
\ \ \lambda\in(1,+\infty)$$ on $\mathbb{C}^n$ which are expanding
K$\ddot{a}$hler-Ricci solitons, where $u_\lambda:\
\mathbb{R}\rightarrow\mathbb{R}$ is a suitable smooth convex and
increasing function for each $\lambda\in(1,+\infty)$. It was also
shown by Cao [1] that these K$\ddot{a}$hler metrics have positive
sectional curvature everywhere on $\mathbb{C}^n$. It was further
shown in [5] that these K$\ddot{a}$hler metrics have positive
curvature in the sense of Nakano (i.e., positive curvature
operator on the subspace of (1,1)-forms) and the geodesic distance
from the origin of $\mathbb{C}^n$ with respect to the metric
$g^\lambda_{i\bar{j}}$ (for each $\lambda\in(1,+\infty)$) is
asymptotic to $\sqrt{\lambda}|z|^\frac{1}{\lambda}$ as
$|z|\rightarrow+\infty$. Thus for each $\epsilon\in(0,1)$, we
choose properly a $\lambda\in(1,+\infty)$ such that the
K$\ddot{a}$hler manifold $M^n=(\mathbb{C}^n,g^\lambda_{i\bar{j}})$
acquires the equalities for all $d$ in the improved dimension
estimate of Theorem 3.1. \vskip 1cm

\section{Generalizations}\indent

The main purpose of this section is to prove Theorem 1.4 in the
introduction. Let us begin with a lemma on exhausting function
which does not demand a upper bound on the curvature. \vskip
0.5cm\noindent {\bf Lemma 4.1}{\sl \ \ Let $M^n$ be a complete
noncompact K$\ddot{a}$hler manifold of complex dimension $n$ with
nonnegative holomorphic bisectional curvature. Then there exists a
positive constant $C(n)$ depending only on the dimension $n$ such
that for any $a>0$, there is a smooth function $\varphi$ on $M^n$
satisfying
$$C(n)^{-1}(1+\frac{r(x,x_0)}{a})\leq\varphi(x)\leq
C(n)(1+\frac{r(x,x_0)}{a}),\ \ \ x\in M^n,$$
$$|\nabla\varphi|\leq\frac{C(n)}{a},\ \ \ on\ M^n,$$
$$|\varphi_{i\bar{j}}|\leq\frac{C(n)}{a^2},\ \ \ on\ M^n,$$
where $x_0$ is a fixed point in $M^n$ and $r(x,x_0)$ is the
geodesic distance between $x$ and $x_0$.} \vskip 0.1cm\noindent
{\sc Proof.}\ \ \ Consider the following heat equation
$$\left\{\begin{array}{l}
\frac{\partial u}{\partial t}=\Delta u,\ \ \ on\
M^n\times(0,+\infty),\\[3mm]
u|_{t=0}=r(x,x_0)+1,\ \ \ on\ M^n.
\end{array}\right.\eqno(4.1)$$
It is clear (see for example Lemma 1.2 in [15]) that the solution
$u(x,t)$ exists for all $t\in(0,+\infty)$ and $u(x,t)$ can be
represented by
$$u(x,t)=\int_{M^n}H(x,y,t)(r(y,x_0)+1)dy.$$
By the heat kernel estimate of Li-Yau [13], it is not hard (see
for example Corollary 1.4 in [15]) to see
$$C(n)^{-1}(1+r(x,x_0))\leq u(x,1)\leq
C(n)(1+r(x,x_0)),\ \ \ for\ x\in M^n,\eqno(4.2)$$ for some
positive constant $C(n)$ depending only on the dimension $n$.

From the heat equation (4.1) and the complex
Lichnerowicz-Laplacian heat equation (2.8), it is not hard to see
$$
\arraycolsep=1.5pt
\begin{array}{rcl}
(\frac{\partial}{\partial t}-\Delta)|\nabla
u|^2&=&-|u_{i\bar{j}}|^2-|u_{ij}|^2-2R_{i\bar{j}}u^iu^{\bar{j}}\\[3mm]
&\leq&-|u_{i\bar{j}}|^2,
\end{array}\eqno(4.3)
$$
$$$$
$$
\arraycolsep=1.5pt
\begin{array}{rcl}
(\frac{\partial}{\partial t}-\Delta)|u_{i\bar{j}}|^2&=&-(|u_{i\bar{j},k}|^2+|u_{i\bar{j},\bar{k}}|^2)
                  -R_{i\bar{i}j\bar{j}}(\lambda_i-\lambda_j)^2\\[3mm]
&\leq&0,
\end{array}\eqno(4.4)
$$
where $\lambda_i$, $i=1,\cdots,n$, are the eigenvalues of
$u_{i\bar{j}}$ with respect to the metric $g_{i\bar{j}}$ and we
have used the curvature condition in above two inequalities.

The maximum principles for the solutions of the heat equation
(4.1) with at most exponential growth on complete noncompact
K$\ddot{a}$hler manifolds with nonnegative holomorphic bisectional
curvature have been well developed by Ni-Tam in [15] (especially,
Theorem 2.1 and Theorem 3.1 and their proofs in [15]). In
particular, the maximum principles for $u(x,t)$ and $u_{i\bar{j}}$
work in our present case because the solution $u(x,t)$ is at most
linear growth. Since $|\nabla r|=1$, by the equation (4.3) we have
$$|\nabla u|(x,t)\leq1,\ \ \ for\ all\ x\in M^n,\ t>0.\eqno(4.5)$$
To estimate $|u_{i\bar{j}}|$, let
$v(x,t)=t|u_{i\bar{j}}(x,t)|^2+|\nabla u(x,t)|^2$. By a direct
computation, we have
$$\left\{\begin{array}{l}
(\frac{\partial}{\partial t}-\Delta)v(x,t)\leq0,\ \ \ on\
M^n\times(0,+\infty),\\[3mm]
v|_{t=0}\leq1,\ \ \ on\ M^n.
\end{array}\right.$$
which implies that
$$|u_{i\bar{j}}|^2(x,1)\leq1,\ \ \ for\ all\ x\in M^n.\eqno(4.6)$$

Summarizing (4.2), (4.5) and (4.6), we obtain the desired function
$\varphi(x)=u(x,1)$ for $a=1$. For general $a>0$, we consider the
metric $\tilde{g}_{i\bar{j}}=\frac{1}{a^2}g_{i\bar{j}}$ on $M^n$.
Obviously the holomorphic bisectional curvature is still
nonnegative. By applying the result for the case $a=1$, there
exists a function $\varphi_a$ on $(M^n,\tilde{g}_{i\bar{j}})$ such
that
$$C(n)^{-1}(1+\tilde{r}(x,x_0))\leq\varphi_a(x)\leq
C(n)(1+\tilde{r}(x,x_0)),$$
$$|\tilde{\nabla}\varphi_a|_{\tilde{g}_{i\bar{j}}}\leq C(n),$$
$$|(\varphi_a)_{i\bar{j}}|_{\tilde{g}_{i\bar{j}}}\leq C(n),$$
on $M^n$. Since $\tilde{r}(x,x_0)=\frac{1}{a}r(x,x_0)$,
$|\tilde{\nabla}\varphi_a|_{\tilde{g}_{i\bar{j}}}=a|\nabla\varphi_a|_{g_{i\bar{j}}}$
and
$|(\varphi_a)_{i\bar{j}}|_{\tilde{g}_{i\bar{j}}}=a^2|(\varphi_a)_{i\bar{j}}|_{g_{i\bar{j}}}$,
the function $\varphi_a(x)$ fulfills all the requirements of the
lemma. Therefore the proof of Lemma 4.1 is completed.\hskip 1cm\#
\vskip 0.5cm It is clear that for the exhausting function
$\varphi$ in Lemma 4.1, there holds
$$B(x_0,\frac{a}{C(n)})\subset\{\varphi\leq C(n)+1\}\subset
B(x_0,C(n)(C(n)+1)a).\eqno(4.7)$$ We are now ready to prove
Theorem 1.4.\\
{\bf Proof of Theorem 1.4}

Let $M^n$ be a complete noncompact K$\ddot{a}$hler manifold of
complex dimension $n$ with nonnegative holomorphic bisectional
curvature satisfying
$$Vol(B(x_0,r))\leq C(1+r)^{2k},\ \ \ for\ all\
r\geq0,\eqno(4.8)$$ for some $1\leq k\leq n$. In views of Theorem
1.2, we may assume $k<n$. Let $\mathcal{P}$ denote the set of all
smooth plurisubharmonic functions of at most logarithmic growth.
For the fixed point $x_0\in M^n$, define a number
$$r_{x_0}=\max\{r|\ h\in\mathcal{P},\ (h_{i\bar{j}}(x_0))\ has\ rank\
r\}.$$ Obviously $0\leq r_{x_0}\leq n$. We first claim
$$r_{x_0}\leq[k]<n.\eqno(4.9)$$
Suppose not, then there exists a plurisubharmonic function
$h\in\mathcal{P}$ such that $(h_{i\bar{j}}(x_0))$ has rank at
least $[k]+1$. Obviously by adding a constant, we may assume
$h(x_0)\geq2$. We modify the function $h$ as follows. First we set
$$ \tilde{h}=\max\{h,\ 1\}\geq1$$
which is a continuous plurisubharmonic function of logarithmic
growth. Next we evolve the function $\tilde{h}$ by the heat
equation to get a solution $h(x,t)$ on $M^n\times(0,+\infty)$. By
applying the results of Ni-Tam [15] (see Theorem 3.1, Theorem 2.1
and Corollary 1.4 in [15]) we know that for fixed $t=1$, the
function $h(x,1)$ is a smooth plurisubharmonic function on $M^n$,
$(h_{i\bar{j}}(x_0,1))$ has also rank $r_{x_0}\geq[k]+1$, and
$h(x,1)$ is still of logarithmic growth and satisfies
$$1\leq h(x,1)\leq Clog(2+r(x,x_0)),\ \ \ on\ M^n,\eqno(4.10)$$
for some positive constant $C$. In the following we use this new
function $h(x,1)$ to replace the original one and still denote it
by $h$.

Since $(h_{i\bar{j}}(x_0))$ has rank at least $[k]+1$, we have
$$(\sqrt{-1}\partial\bar{\partial}h)^{[k]+1}\wedge\omega^{n-([k]+1)}>0,\
\ \ in\ a\ small\ neighborhood\ of\ x_0\in M^n,$$ where
$\omega=\sqrt{-1}g_{i\bar{j}}dz^i\wedge d\bar{z}^j$ is the
K$\ddot{a}$hler form associated to the K$\ddot{a}$hler metric
$g_{i\bar{j}}$. For any $a>2$ and let $\varphi$ be the exhausting
function obtained in lemma 4.1. Then there exists a positive
constant $\delta>0$ independent of $a$ such that
$$0<\delta\leq\int_{\{\varphi\leq
C(n)+1\}}(1-\frac{\varphi}{(C(n)+1)})^{2([k]+1)}(\sqrt{-1}\partial\bar{\partial}h)^{[k]+1}
\wedge\omega^{n-([k]+1)}\eqno(4.11)$$ On the other hand, by Lemma
4.1, (4.7) and (4.10), we get from integrating by parts
\begin{displaymath}
\begin{split}
   &\int_{\{\varphi\leq C(n)+1\}}(1-\frac{\varphi}{(C(n)+1)})^{2([k]+1)}
    (\sqrt{-1}\partial\bar{\partial}h)^{[k]+1}\wedge\omega^{n-([k]+1)}\\[4mm]
  \leq &\frac{\tilde{C}(n)}{a^2}\int_{\{\varphi\leq
    C(n)+1\}}(1-\frac{\varphi}{(C(n)+1)})^{2[k]}h
    (\sqrt{-1}\partial\bar{\partial}h)^{[k]}\wedge\omega^{n-[k]}\\[4mm]
\leq &\frac{C\tilde{C}(n)loga}{a^2}\int_{\{\varphi\leq
    C(n)+1\}}(1-\frac{\varphi}{(C(n)+1)})^{2[k]}
    (\sqrt{-1}\partial\bar{\partial}h)^{[k]}\wedge\omega^{n-[k]}\\[4mm]
\leq
 &\cdots\leq\frac{C^{([k]+1)}\tilde{C}(n)(loga)^{([k]+1)}}{a^{2([k]+1)}}\int_{\{\varphi\leq
    C(n)+1\}}\omega^n\\[4mm]
\leq&\frac{C^{([k]+1)}\tilde{C}(n)(loga)^{([k]+1)}}{a^{2([k]+1)}}Vol(B(x_0,C(n)(C(n)+1)a)),
\end{split}
\end{displaymath}
$$\eqno(4.12)$$
where $\tilde{C}(n)$ denotes various positive constants depending
only on the dimension $n$.

The combination of (4.11) and (4.12) gives
$$Vol(B(x_0,a))\geq C\frac{a^{2([k]+1)}}{(loga)^{([k]+1)}},\ \ \ \forall\ a\geq2,$$
for some positive constant $C$ independent of $a$. Since
$[k]+1>k$, this arrives a contradition with the volume growth
assumption (4.8). Thus we have proved the claimation (4.9).

Let $h$ be a smooth plurisubharmonic function of logarithmic
growth such that $(h_{i\bar{j}}(x_0))$ has the maximal rank
$r_{x_0}$. Evolve the function $h$ by the heat equation
$$\left\{\begin{array}{l}
\frac{\partial u}{\partial t}=\Delta u,\ \ \ on\
M^n\times(0,+\infty)\\[3mm]
u|_{t=0}=h,\ \ \ on\ M^n.
\end{array}\right.$$
Clearly the smooth solution $u(x,t)$ exists for all
$t\in(0,+\infty)$. And by applying the results of Ni-Tam [15]
again, we know that for each fixed $t>0$, the function
$u(\cdot,t)$ is still a smooth plurisubharmonic function of
logarithmic growth on $M^n$ and $(u_{i\bar{j}}(x_0,t))$ still has
the maximal rank $r_{x_0}$. Then as before by using Corollary 2.1
in [15], we know that the kernel space $K(x,t)$ of
$(u_{i\bar{j}}(x,t))$ is a distribution which is invariant under
parallel translations and then $(u_{i\bar{j}}(x,t))$ has rank
$r_{x_0}$ everywhere, moreover the universal cover $\tilde{M}^n$
splits isometrically and holomorphically as
$$\tilde{M}^n=M_1^{n-r_{x_0}}\times M_2^{r_{x_0}}\eqno(4.13)$$
where $K$ corresponds the tangent bundle of $M_1^{n-r_{x_0}}$ and
$(u_{i\bar{j}}(x,t))>0$ on $M_2^{r_{x_0}}\times(0,+\infty)$. Both
$M_1^{n-r_{x_0}}$ and $M_2^{r_{x_0}}$ are complete K$\ddot{a}$hler
manifold with nonnegative holomorphic bisectional curvature. The
estimate (4.9) shows that the complex dimension of the kernel $K$
is at least $n-[k]\geq1$.

Clearly for any $f\in\mathcal{O}_d(M^n)$, we can lift $f$ as a
function in $\mathcal{O}_d(\tilde{M}^n)$, still denoted by $f$. In
the following we will show that such $f$ must descend to a
function of $\mathcal{O}_d(M_2^{r_{x_0}})$.

For any nonconstant $f\in\mathcal{O}_d(M^n)$, let $u_f(x,t)$ solve
the heat equation
$$\left\{\begin{array}{l}
\frac{\partial}{\partial t}u_f=\Delta u_f,\ \ \ on\
M^n\times(0,+\infty)\\[3mm]
u_f|_{t=0}=log(|f|^2+1),\ \ \ on\ M^n.
\end{array}\right.$$
We see as before by applying the results of Ni-Tam [15] that the
solution $u_f$ is smooth plurisubharmonic with logarithmic growth
on $M^n$ for each $t\in(0,+\infty)$ and there exists $\tilde{t}>0$
such that for any $0<t<\tilde{t}$, the kernel space $K_f(x,t)$ is
a distribution which is invariant under parallel translations and
the universal cover $\tilde{M}^n$ splits isometrically and
holomorphically as
$$\tilde{M}^n=(M_f)_1^{n-q}\times(M_f)_2^q$$
with $((u_f)_{i\bar{j}})|_{(M_f)_1^{n-q}}\equiv0$ and
$((u_f)_{i\bar{j}})|_{(M_f)_2^q}>0$. Here we have lifted the
function $u_f$ to the universal cover and still denoted it by
$u_f$. Recall that $(u_{i\bar{j}}(x,t))$ has maximal rank
$r_{x_0}$ everywhere. It follows that the orthogonal completement
$K_f^\perp(x_0,t)$ of $K_f(x_0,t)$ is contained in the orthogonal
completement $K^\perp(x_0,t)$ of $K(x_0,t)$, otherwise
$(u+u_f)_{i\bar{j}}(x_0,t)\ (t\in(0,\tilde{t}))$ would have rank
$\geq$ $r_{x_0}+1$, which contradicts with the definition of
$r_{x_0}$. And then by the parallel translation invariance of
$K_f^\perp$ and $K^\perp$ we have
$$K_f^\perp\subset K^\perp$$
and
$$(M_f)_2^q\subset M_2^{r_{x_0}}.\eqno(4.14)$$
Let us write $x=(x_1,x_2)\in M^n=(M_f)_1^{n-q}\times(M_f)_2^q$
with $x\in(M_f)_1^{n-q}$ and $x_2\in(M_f)_2^q$. Since
$((u_f)_{i\bar{j}})|_{(M_f)_1^{n-q}}\equiv0$, we know that for
each fixed $t\in(0,\tilde{t})$ and arbitrarily fixed
$x_2\in(M_f)_2^q$, the function $u_f(\cdot,x_2,t)$ is a harmonic
function of logarithmic growth. We then conclude from Cheng-Yau
[4] that for each fixed $t\in(0,\tilde{t})$, $u_f$ is only a
function of the second variable $x_2\in(M_f)_2^q$. This implies
from Yau's Liouville theorem [19] that $f$ is only a function of
variable $x_2\in(M_f)_2^q$. Thus by combining with (4.14), we have
$$\mathcal{O}_d(M^n)\subset\mathcal{O}_d(M_2^{r_{x_0}}),\eqno(4.15)$$
in particular,
$$dim_{\mathbb{C}}\mathcal{O}_d(M^n)\leq
dim_{\mathbb{C}}\mathcal{O}_d(M_2^{r_{x_0}}).\eqno(4.16)$$
Therefore by using Theorem 1.2 and (4.9), we deduce that
$$dim_{\mathbb{C}}\mathcal{O}_d(M^n)\leq
dim_{\mathbb{C}}\mathcal{O}_{[d]}(\mathbb{C}^{[k]}),\ \ \ for\
each\ d>0.\eqno(4.17)$$

We next discuss the rigidity part of Theorem 1.4. Suppose there
exists a positive integer $d$ such that
$$dim_{\mathbb{C}}\mathcal{O}_d(M^n)=dim_{\mathbb{C}}\mathcal{O}_d(\mathbb{C}^{[k]}).$$
From (4.16) and (4.17) we know that
$$r_{x_0}=[k]$$
and
$$dim_{\mathbb{C}}\mathcal{O}_d(M_2^{r_{x_0}})=dim_{\mathbb{C}}\mathcal{O}_d(\mathbb{C}^{[k]}).$$
Thus by Theorem 1.2 we deduce that $M_2^{r_{x_0}}$ is
holomorphically isometric to $\mathbb{C}^{[k]}$.\\
So
$$\tilde{M}^n=M_1^{n-[k]}\times \mathbb{C}^{[k]}\eqno(4.18)$$
isometrically and holomorphically. And by combining with the
inclusion (4.15) we have
$$\mathcal{O}_d(M^n)=\mathcal{O}_d(\tilde{M}^n)=\mathcal{O}_d(\mathbb{C}^{[k]}).$$
This says that every function in $\mathcal{O}_d(\mathbb{C}^{[k]})\
(\subset\mathcal{O}_d(\tilde{M}^n))$ is $\pi_1(M^n)$-invariant. In
particular, the coordinate functions $z^1,\cdots,z^{[k]}$ of
$\mathbb{C}^{[k]}$ (regarding as holomorphic functions of
$\tilde{M}^n=M_1^{n-[k]}\times\mathbb{C}^{[k]}$) are
$\pi_1(M^n)$-invariant.

Let $\sigma\in\pi_1(M^n)$ be any deck transformation of
$\tilde{M}^n$ and let $x_2^0=(z_0^1,\cdots,z_0^{[k]})$ be any
fixed point in $\mathbb{C}^{[k]}$. By the $\pi_1(M^n)$-invariance
of $z^1,\cdots,z^{[k]}$, we have
$$\sigma(\bigcap\limits_{i=1}^{[k]}[z^i=z_0^i])=\bigcap\limits_{i=1}^{[k]}[z^i=z^i_0],$$
i.e.,
$$\sigma(M_1^{n-[k]}\times\{x_2^0\})=M_1^{n-[k]}\times\{x_2^0\}.$$
Then for each fixed $x_2^0\in\mathbb{C}^{[k]}$, the deck
transformation induces a transformation, denoted by
$\sigma_{x_2^0}$, on $M_1^{n-[k]}$. Let $x_1\in M_1^{n-[k]}$, and
$x_2$, $x_2'\in\mathbb{C}^{[k]}$ be arbitrary. Then the reduced
transformation $\sigma_{x_2}$ and $\sigma_{x_2'}$ are given by
$$\sigma(x_1,x_2)=(\sigma_{x_2}(x_1),x_2),$$
and
$$\sigma(x_1,x_2')=(\sigma_{x_2'}(x_1),x_2').$$
Since the deck transformation $\sigma$ is isometric on
$\tilde{M}^n=M_1^{n-[k]}\times\mathbb{C}^{[k]}$, we deduce that
$$\sigma_{x_2}(x_1)=\sigma_{x_2'}(x_1),\ \ \ for\ all\ x_1\in M_1^{n-[k]}.$$
This implies that the action $\pi_1(M^n)$ on
$\tilde{M}^n=M_1^{n-[k]}\times\mathbb{C}^{[k]}$ splits and acts
trivially on the second factor $\mathbb{C}^{[k]}$. Hence
$$M^n=(M_1^{n-[k]}/\pi_1(M^n))\times\mathbb{C}^{[k]}$$
isometrically and holomorphically. Clearly the dimension estimate
(4.17) implies that the first factor $(M_1^{n-[k]}/\pi_1(M^n))$
can not carry any nonconstant holomorphic function of polynomial
growth.

Therefore we have completed the proof of Theorem 1.4.\hskip 1cm\#
\vskip 0.5cm Finally it is not hard to see that the combination of
the arguments of Section 3 and the above proof of Theorem 1.4
gives the following improved dimension estimate.\\
{\bf Proposition 4.2} \ \ {\sl Let $M^n$ be a complete noncompact
K$\ddot{a}$hler manifold of complex dimension $n$ with nonnegative
holomorphic bisectional curvature. Assume that its Ricci curvature
is positive at least at one point in $M^n$ and suppose there exist
a point $x_0\in M$, and positive constants $1\leq k\leq n$ and
$C>0$ such that
$$Vol(B(x_0,r))\leq C(1+r)^{2k},\ \ \ for\ all\ r\geq0.$$
Then there exists a positive constant $\epsilon\in(0,1)$ such that
$$dim_{\mathbb{C}}\mathcal{O}_d(M^n)\leq
dim_{\mathbb{C}}\mathcal{O}_{[(1-\epsilon)d]}(\mathbb{C}^{[k]}),$$
for all positive integers $d$.\hskip 1cm\#

\end{document}